\documentclass[12pt]{article}
\usepackage{amsmath,amssymb,amsfonts,graphicx,wrapfig}

\newtheorem{theorem}{Theorem}
\newtheorem{lemma}{Lemma}

\begin{document}

{\large
\begin{center}
Existential monadic second order convergence law fails on sparse random graphs
\end{center}
}

\begin{center}
Alena Egorova, Maksim Zhukovskii\footnote{Moscow Institute of Physics and Technology, laboratory of advanced combinatorics and network applications}
\end{center}

\kern 1.5 cm

\begin{center}
{\bf Abstract}
\end{center}

In the paper, we prove that existential monadic second order convergence law fails for the binomial random graph $G(n,n^{-\alpha})$ for every $\alpha\in(0,1)$.

\kern 1.5 cm

\section{Introduction}
\label{s1}

Sentences in the first order language of graphs (FO sentences) are constructed using relational symbols $\sim$ (interpreted as adjacency) and $=$, logical connectives $\neg,\rightarrow,\leftrightarrow,\vee,\wedge$, variables $x,y,x_1, \ldots$ that express vertices of a graph, quantifiers $\forall,\exists$ and parentheses.  Monadic second order, or MSO, sentences are built of the above symbols of the first order language, as well as the variables $X,Y,X_1,\ldots$ that are interpreted as unary predicates. In an MSO sentence, variables $x,y,x_1, \ldots$ (that express vertices) are called {\it FO variables}, and variables $X,Y,X_1,\ldots$ (that express sets) are called {\it MSO variables}. If, in an MSO sentence $\phi$, all the MSO variables are existential and in the beginning (that is 
\begin{equation}
\phi=\exists X_1\ldots\exists X_m\,\,\varphi(X_1,\ldots,X_m)
\label{general_EMSO}
\end{equation} 
where $\varphi(X_1,\ldots,X_m)$ is a FO sentence with unary predicates $X_1,\ldots,X_m$), then the sentence is called {\it existential} monadic second order (EMSO). Sentences must have finite number of logical connectivities. In what follows, for a sentence $\phi$, we use the usual notation from model theory $G\models\phi$ if $\phi$ is true for $G$. More detailed (and more formal) definitions can be found, e.g., in \cite{Libkin,Survey}.\\

In this paper, we consider the binomial model of random graph $G(n,p)$. In this model, we have $G(n,p) = (V_n,E)$, where $V_n=\{1,\ldots,n\}$, and each pair of vertices is connected by an edge with probability $p$ and independently of other pairs. For more information, we refer readers to the books \cite{AS,Bollobas,Janson}. Clearly, $G(n,\frac{1}{2})$ is distributed uniformly on the set of all graphs on $V_n$. Y. Glebskii, D. Kogan, M. Liogon'kii and V. Talanov in 1969~\cite{Glebskii}, and independently R. Fagin in 1976~\cite{Fagin}, proved that any FO sentence is either true with asymptotical probability 1 (asymptotically almost surely or a.a.s.) or a.a.s. false for $G(n,\frac{1}{2})$, as $n\to \infty$, i.e. $G(n,\frac{1}{2})$ \textit{obeys the FO zero-one law}. For MSO, the zero-one law was disproved by M. Kaufmann and S. Shelah in 1985~\cite{Kaufmann_Shelah}. They prove that there is even no {\it MSO convergence law} (i.e., there is an MSO sentence $\phi$ such that ${\sf P}(G(n,\frac{1}{2})\models\phi)$ does not converge). After that, in 1987~\cite{Kaufmann}, Kaufmann proved that there exists an EMSO sentence with {\it 4 binary relations} ({\it undirected graphs} that are considered in this paper have only one symmetric binary relation apart from the default relation $=$) that has no asymptotic probability. The non-convergence result for $G(n,\frac{1}{2})$ was obtained by J.-M. Le Bars in 2001~\cite{Le_Bars}. 

The binomial random graph $G(n,n^{-\alpha})$ (where $\alpha\in(0,1)$ is a fixed constant) is called {\it sparse}. In 1988~\cite{Shelah}, S. Shelah and J. Spencer studied FO logic of sparse random graphs and proved that FO zero-one law holds if and only if $\alpha$ is irrational. For MSO, for all $\alpha\in(0,1)$, the 0-1 law (and even the convergence law) was disproved by J. Tyszkiewicz in 1993~\cite{Tyszk}. EMSO sentences are considered in the very brief last chapter of this paper (where, in particular, Tyszkiewicz mentioned that, for $G(n,1/2)$, the EMSO convergence law was disproved by Kaufmann in~\cite{Kaufmann}, but this is false). Moreover, he claimed that the techniques from the previous chapters of~\cite{Tyszk} can be also applied for disproving EMSO convergence, but did not give the proof. He only mentioned that the proof is based on a certain property of the random graph, which can be verified using some well known combinatorial estimates. We have tried to check this idea, but we failed. In the present paper, we use another method of constructing sentences that have non-convergent probabilities and prove the non-convergence for all $\alpha$ in the range $(0,1)$. One of the advantages of our method is that, for $\alpha<\frac{1}{2}$, it gives an EMSO sentence with only one monadic variable.\\

Let us finish this section with the following remark. We do not consider $\alpha\geq 1$ (which is usually referred as the {\it very sparse} case) in this paper, since this case is completely closed. In~\cite{Shelah} it is proven, that $G(n,n^{-1-1/m})$, for positive integers $m$, does not obey FO zero-one law (nevertheless, FO convergence law holds in this case~\cite{Lynch}). For all the remaining $\alpha>1$ FO zero-one holds. Finally, the random graph $G(n,1/n)$ does not obey FO zero-one law, but there is convergence~\cite{Lynch}. The same results are true for MSO logic~\cite{Luczak,Zhuk_Ostr_Annals} as well. Since FO$\subset$EMSO$\subset$MSO, the same is also true for EMSO logic.

\section{The non-convergence result}
\label{result}

The main result of this paper is given below.

\begin{theorem}
Let $\delta>0$. 
\begin{enumerate}
\item There exists an EMSO sentence $\phi$ with 1 monadic variable such that, for every $\alpha<\frac{1}{2}-\delta$, ${\sf P}(G(n,n^{-\alpha}) \models\phi)$ does not converge as $n \to \infty$. 
\item There exists an EMSO sentence $\phi$ such that, for every $\alpha<1-\delta$, ${\sf P}(G(n,n^{-\alpha}) \models\phi)$ does not converge as $n \to \infty$. 
\end{enumerate}
\label{main}
\end{theorem}

In particular, Theorem~\ref{main} states that, for every $\alpha\in(0,1)$, EMSO convergence law fails for $G(n,n^{-\alpha})$. The proof of this result is constructive, and, for $\alpha<\frac{1}{2}$, the obtained construction is very short in terms of the number of monadic variables. It is also worth mentioning that, given an arbitrary small $\varepsilon>0$, one construction works for all $\alpha<1-\varepsilon$.

Note that the constructions of Kaufmann~\cite{Kaufmann}, Le Bars~\cite{Le_Bars} and Tyszkiewicz~\cite{Tyszk} (the latter is not existential) exploit more monadic variables. In particular, the construction of Kaufmann has 4 monadic variables, and the modified construction of Le Bars has even more. The approach of Tyszkiewicz in~\cite{Tyszk}  is very different, and it requires much more variables. He does not give an explicit construction, and we have not tried to find an optimal way of expressing it. Nevertheless, even the explicit part of this construction contains 7 monadic variables.\\

The scheme of the proof of Theorem~\ref{main} is the following. Let $p=n^{-\alpha}$. For $\alpha<\frac{1}{2}$, we consider a certain graph sequence $H_j$, $j\in\mathbb{N}$, such that the number of vertices $v(H_j)\sim C4^j$ (constant before the exponent does not play a role in the proof of non-convergence --- it appears because of the further condition on the graph sequence) and `being isomorphic to $H_j$ for some $j$' is FO expressible. If so, we may construct a sentence $\phi=\exists X\,\varphi(X)\wedge\max(X)$, where the FO formula $\varphi(X)$ says that, for some $j$, a subgraph induced on $[X]:=\{v:\,X(v)\}$ is isomorphic to $H_j$, and the formula $\max(X)$ says that every vertex outside $[X]$ has a neighbor inside $[X]$. 

After that, we prove that our graph sequence is so nice that there exists a $k=k(n)$ and $\kappa\in(1/2,1)$ with the following three properties: 1) if $v(H_j)>k+\varepsilon$, then it is not likely that $G(n,p)$ contains an induced $H_j$; 2) if  $v<(\kappa-\varepsilon)k$, then it is not likely that $G(n,p)$ contains a subset of size at most $v$ such that every vertex outside this set has a neighbor inside; and, finally, 3) consider an infinite sequence of $n$ such that there exists $j$ with $v(H_j)\in((\kappa+\varepsilon)k,(1-\varepsilon)k)$, then it is likely that, for these $n$, $G(n,p)$ contains an induced $H_j$ such that every vertex outside this subgraph has a neighbor inside.

The important phenomenon which plays a key role in our proof (and which is described in properties 1) and 2)) is that a size of a maximum induced $H_j$ and a size of a minimum set which has no `outside-isolated' vertices are close to each other but differs in $\kappa\in(1/2,1)$ factor (the second extremum is smaller). Therefore, if, for infinitely many $n$, there are no $v(H_j)$ between these thresholds, this give us the result. And the latter is true since, for large $j$, $v(H_{j+1})/v(H_j)$ equals roughly $4$, while $\kappa^{-1}<2$. The detailed proof of the first part of Theorem~\ref{main} is given in Section~\ref{before_one_half}.

Unfortunately, for $\alpha\geq\frac{1}{2}$, the same technique does not work because the variance of the number of induced $H_j$ becomes excessively large (it makes it impossible to apply Chebyshev's inequality for the property 3)). It is also impossible to apply martingale techniques (in contrast, it works in a similar situation~\cite{Frieze}) since changing the edges incident with a single vertex may destroy a major part of $H_j$. Luckily, we find a modified graph sequence $H_j^*$ such that $v(H_j^*)\sim C r^{\Theta(r^j)}$ for an appropriate choice of an integer parameter $r$. This modification makes it possible to improve strongly the difference between the mentioned thresholds. In Section~\ref{after_one_half}, we define this sequence and give the proof of the second part of Theorem~\ref{main}.

\subsection{Non-converegence for $\alpha<1/2$}
\label{before_one_half}

\subsubsection{The graph sequence}

\begin{figure}[h!]
\center{\includegraphics[scale=0.2]{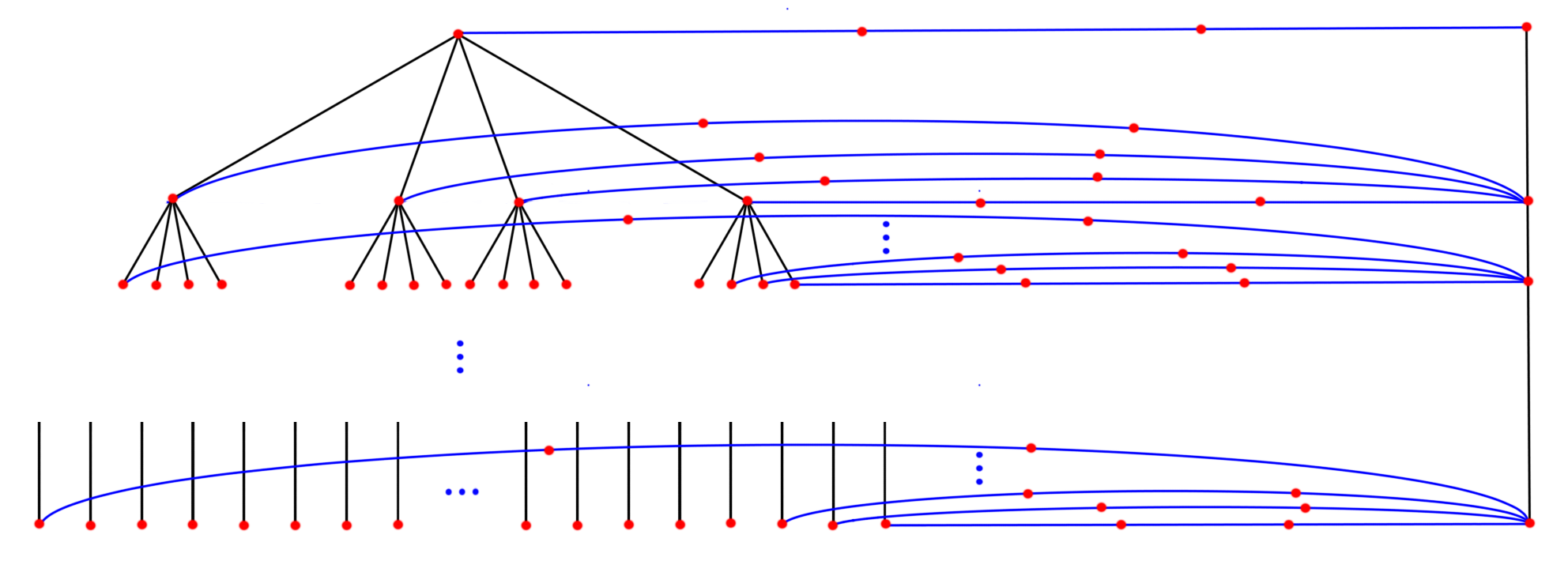}}
\caption{2-product of a simple path and a perfect $4$-ary tree.}
\label{fig:the_graph}
\end{figure}

Consider two rooted trees $F_1$ and $F_2$ with roots $R_1$ and $R_2$ respectively and a non-negative integer $\gamma$. Let us define {\it the $\gamma$-product} of the rooted trees $F_1\cdot^{\gamma} F_2$ in the following way. Let $\mathcal{E}$ be the set of all pairs $(u,v)$ where $u\in V(F_1)$, $v\in V(F_2)$ and $u$, $v$ are at the same distance from $R_1$, $R_2$ in $F_1$, $F_2$ respectively. Then $F_1\cdot^{\gamma} F_2$ is obtained from the disjoint union $F_1\sqcup F_2$ in the following way: for every pair $(u,v)\in\mathcal{E}$, we add to the graph a $P_{\gamma+2}$ (i.e., simple path with $\gamma+1$ edges) connecting $u$ and $v$. Fix an arbitrary positive integer $a$ and consider two trees $F_1,F_2$ with $a$ and $\frac{4^a-1}{3}$ vertices respectively: $F_1$ is a simple path rooted at one of its end-points, and $F_2$ is a perfect $4$-ary tree (every non-leaf vertex of $F_2$ has 4 children and, for every $i\in\{1,\ldots,a-1\}$, the number of vertices at the distance $i$ from $R$ equals $4^i$) rooted at the only vertex having degree $4$. Denote $W_a^{\gamma}=F_1\cdot^{\gamma} F_2$ (see Fig.~\ref{fig:the_graph}).

Surely, $W_a^{\gamma}$ has $a+(\gamma+1)\frac{4^a-1}{3}$ vertices and $a+(\gamma+2)\frac{4^a-1}{3}-2$ edges.

\subsubsection{The sentence}
\label{sentence_desc}

In~\cite{Zhuk_Le_Bars} (see the proof of  Theorem 3), we construct a FO formula $\varphi_0(X)$ with two binary predicates $\sim,=$ and one unary predicate $X$ saying that `for some $a$, the induced subgraph on $[X]$ is isomorphic to $W_a^0$'. In the same way, it is straightforward to construct a FO formula $\varphi_{\gamma}(X)$ saying the same but about $W_a^{\gamma}$. For reader's convenience, let us briefly describe this sentence:
$$
 \varphi_{\gamma}=\exists x\exists y^1\exists y_1^2\ldots\exists y_4^2\exists z_1\exists z_2\exists w\exists h\quad
 [\mathrm{DEG}(x,\ldots,h)\wedge\mathrm{PATH}(z_2,w)\wedge\mathrm{PROD}(y_1,y^2_1,\ldots,y^2_4)\wedge
$$
$$
 \mathrm{PERFECT}(z_1,z_2,y_1^2,\ldots,y_4^2)\wedge\mathrm{TREE}(x,z_1)\wedge\mathrm{START}(z_1,z_2)
 \wedge\mathrm{MAX}],
$$
where $x$ is the end-point (root) of the simple path $F_1$ which is connected via $P_{\gamma+2}$ to the root $y^1$ of $F_2$ (in this case, we say that $x$ is $\gamma$-neighbor of $y^1$); $y_1^2,\ldots,y_4^2$ are children of $y^1$ in $F_2$; $z_1$ is the child of $x$ in $F_1$; $z_2$ is the child of $z_1$ in $F_1$; $w$ is the second end-point of $F_1$, and $h$ is its parent.

Here,

\begin{itemize}

\item $\mathrm{DEG}(x,\ldots,h)$ 

--- says that the vertex set of $[X]$ can be partitioned in 3 parts (vertices of $F_1$, vertices of $F_2$ and inner vertices of $P_{\gamma+2}$'s) by defining degrees of all the vertices in the induced subgraph on the set $[X]$, 

--- defines the relations between the distinguished vertices (in particular, it says that there is a path $P_{\gamma+2}$ between $x$ and $y^1$ such that all its inner vertices have degree $2$ in $[X]$);

\item $\mathrm{START}(z_1,z_2)$ says that each of the four children of the root in $F_2$ has four children in $F_2$, and each of them is connected with $z_2$ via $P_{\gamma+2}$ such that all its inner vertices have degree 2 in $[X]$;

\item $\mathrm{PATH}(z_2,w)$ says that $F_1$ is either a simple path or a union of a simple path and simple cycles;

\item $\mathrm{PROD}(y_1,y^2_1,\ldots,y^2_4)$ says that every vertex $u$ of $F_2$ has an only vertex in $F_1$ which is connected with $u$ via $P_{\gamma+2}$, and all its inner vertices has degree 2 in $[X]$;

\item $\mathrm{PERFECT}(z_1,z_2,y_1^2,\ldots,y_4^2)$ says that every vertex $u$ of $F_2$ (except for the root and the leaves) has a $\gamma$-neighbor in $F_1$ that have two neighbors $v_1$ and $v_2$ in $F_1$ such that $v_1$ has an only $\gamma$-neighbor which is a neighbor of $u$ in $F_2$, and $v_2$ has four $\gamma$-neighbors that are also neighbors of $u$ in $F_2$;

\item $\mathrm{TREE}(x,z_1)$, along with $\mathrm{DEG}(x,\ldots,h)$ and $\mathrm{START}(z_1,z_2)$, says that $F_2$ is a tree;

\item $\mathrm{MAX}$ says that every vertex outside $[X]$ has a neighbor inside $[X]$.

\end{itemize}

For details, see~\cite{Zhuk_Le_Bars}, pages 19--22.\\

Consider an EMSO sentence $\phi_{\gamma}=\exists X\,\varphi_{\gamma}(X)\wedge\max(X)$ with 1 monadic variable, where the FO formula $\max(X)$ with two binary predicates $\sim,=$ and one unary predicate $X$ says that `every vertex outside $[X]$ has a neighbor inside $[X]$'. Then, the result, clearly, follows from the lemma below (the argument is given below after the statement of the lemma).

\begin{lemma}
Let $\alpha\in(0,\frac{1}{2})$. Consider an increasing sequence of positive integers $n_i$. Denote
$$
k_{\gamma}=2\left(1-\frac{\gamma+2}{\gamma+1}\alpha\right),\quad f(x)=x^{\alpha}\ln x.
$$
\begin{enumerate}
\item Let $0<c<\frac{1-\alpha}{k_{\gamma}}$, $\varepsilon>0$. If, for every $i$, there is no integer $a$ such that 
$$
a+(\gamma+1)\frac{4^{a}-1}{3}\in\biggl(c k_{\gamma}f(n_i),k_{\gamma}f(n_i)+\varepsilon\biggr),
$$
then a.a.s., for every $a$, there is no induced copy $F$ of $W_a^{\gamma}$ in $G(n_i,n_i^{-\alpha})$ such that every vertex outside $F$ has a neighbor inside $F$.

\item Let $\gamma>\max\{9,\frac{3\alpha-1}{1-2\alpha}\}$. Let $\frac{1-\alpha}{k_{\gamma}}<C_1<C_2<1$. If, for every $i$, there exists an integer $a_i$ such that $C_1 k_{\gamma}f(n_i)\leq a_i+(\gamma+1)\frac{4^{a_i}-1}{3}\leq C_2 k_{\gamma}f(n_i)$, then a.a.s. there is an induced copy $F$ of $W_{a_i}^{\gamma}$ in $G(n_i,n_i^{-\alpha})$ such that every vertex outside $F$ has a neighbor inside $F$. 
\end{enumerate}
\label{Lem}
\end{lemma}

{\it Remark.} To get the statement of Theorem~\ref{main}.1 from Lemma~\ref{Lem} we need the first part of Lemma~\ref{Lem} to be true only for $\varepsilon=1$. Since its proof does not depend on the value of $\varepsilon$, we state the result in this general form.\\

Now, let us show that Lemma~\ref{Lem} implies Theorem~\ref{main}.1. \\

Let $\frac{1-\alpha}{k_{\gamma}}<C_1<C_2<1$. To prove the result, it is enough to show that there are two sequences $n_i,\,i\in\mathbb{N},$ and $m_i,\,i\in\mathbb{N},$ such that, for all large enough $i$, in 
$[C_1k_{\gamma}f(n_i),C_2 k_{\gamma} f(n_i)]$, there is a number $a_i+(\gamma+1)\frac{4^{a_i}-1}{3}$ for some positive integer $a_i$, and, for $a\in\mathbb{N}$, none of $a+(\gamma+1)\frac{4^{a}-1}{3}$ belongs to $[\frac{4}{5}(1-\alpha)f(m_i),k_{\gamma}f(m_i)+1]$. Indeed, we need to find some positive $c<\frac{1-\alpha}{k_{\gamma}}$ and $\varepsilon>0$ such that, for sufficiently large $i$, there is no integer $a$ where $a+(\gamma+1)\frac{4^{a}-1}{3}\in(c k_{\gamma}f(m_i),k_{\gamma}f(m_i)+\varepsilon)$. Take $c=\frac{4}{5}\frac{1-\alpha}{k_{\gamma}}$ and $\varepsilon=1$.\\

Clearly, $m_i=\lfloor y_i\rfloor$, where $3(1-\alpha) f(y_i)=\frac{1}{3}4^i$, is the desired sequence. This is true since $k_{\gamma}<2(1-\alpha)$ and $(\frac{4}{45}4^i,\frac{2}{9}4^i+1)$ does not contain $a+(\gamma+1)\frac{4^a-1}{3}$ for sufficiently large $i$ and all integers $a$. \\

To find $n_i$, set $C_1<C<C_2$ and take $x_i\in\mathbb{R}$ such that $f(x_i)=\frac{1}{Ck_{\gamma}}\left(i+(\gamma+1)\frac{4^{i}-1}{3}\right)$, $i\in\mathbb{N}$ (note that $\ln x_i\sim \frac{\ln 4}{\alpha}i\to\infty$ as $i\to\infty$).


Setting $n_i:=\lfloor x_i\rfloor$, we get the desired sequence since, for large enough $i$, $Ck_{\gamma}f(x_i)$ belongs to $[C_1k_{\gamma}f(n_i),C_2 k_{\gamma}f(n_i)]$  (indeed, the lower bound follows from $f(n_i)\leq f(x_i)$, and the upper bound follows from $f(n_i)/f(x_i)\geq \frac{f(x_i)-1}{f(x_i)}\to 1$).

\subsubsection{Proof of Lemma~\ref{Lem}}

\noindent 1. First, let us prove that a.a.s., in $G(n,p=n^{-\alpha})$, there is no induced subgraph isomorphic to $W^{\gamma}_a$ for any $a$ such that $a+(\gamma+1)\frac{4^a-1}{3}\geq k_{\gamma}f(n)+\varepsilon$.

Let $a$ be such that $s\geq k_{\gamma}f(n)+\varepsilon$, where $s=a+(\gamma+1)\frac{4^a-1}{3}$.

Let $W(a)$ be the number of induced copies of $W^{\gamma}_a$ in $G(n,p)$. Then 
$$
 {\sf E}W(a)={n\choose s}\frac{s!}{24^{([s-a]/[\gamma+1]-1)/4}}p^{\frac{\gamma+2}{\gamma+1}s-\frac{1}{\gamma+1}a-2}(1-p)^{{s\choose 2}-\frac{\gamma+2}{\gamma+1}s+\frac{1}{\gamma+1}a+2}\leq
$$
$$
 e^{s\left[\left(1-\frac{\gamma+2}{\gamma+1}\alpha\right)\ln n-\frac{sp}{2}-\frac{\ln 24}{4(\gamma+1)}+O(p)\right]+O(\ln^2 n)}.
$$
Therefore, ${\sf P}(\exists a\,\,\text{s.t.}\,\, s\geq k_{\gamma}f(n)+\varepsilon,\,W(a)>0)\leq \sum_{a:\,\,s\geq k_{\gamma}f(n)+\varepsilon}{\sf E}W(a)\to 0$ as $n\to\infty$.\\

It remains to prove that a.a.s., whatever $a$ such that $s\leq ck_{\gamma}f(n)$ is, in $G(n,p)$, there is no induced subgraph $F$ isomorphic to $W^{\gamma}_a$ such that every vertex outside $F$ has a neghbor inside $F$.

We will prove a stronger statement: a.a.s., for every set $X$ of $\lfloor ck_{\gamma}f(n)\rfloor$ vertices, there is a vertex outside $X$ which has no neighbors inside $X$. The probability of this event is at least
$$
 1-{n\choose \lfloor ck_{\gamma}f(n)\rfloor}(1-(1-p)^{\lfloor ck_{\gamma}f(n)\rfloor})^{n-\lfloor ck_{\gamma}f(n)\rfloor}\geq 1-e^{-n^{1-k_{\gamma}c}(1+o(1))+O(n^{\alpha}\ln^2 n)}\to 1
$$
as $n\to\infty$.

\medskip

\noindent 2. Now, let $C_1 k_{\gamma}f(n_i)\leq s_i=a_i+\frac{4^{a_i}-1}{3}\leq C_2 k_{\gamma}f(n_i)$. In what follows, we write $s,a,n$ instead of $s_i,a_i,n_i$ respectively. Let $\tilde W(a)$ be the number of induced copies of $W_a^{\gamma}$ in $G(n,p)$ such that every vertex outside a copy has a neighbor inside. Then
\begin{equation}
 {\sf E}\tilde W(a)={n\choose s}\frac{s!}{24^{([s-a]/[\gamma+1]-1)/4}}p^{\frac{\gamma+2}{\gamma+1}s-\frac{1}{\gamma+1}a-2}(1-p)^{{s\choose 2}-\frac{\gamma+2}{\gamma+1}s+\frac{1}{\gamma+1}a+2}(1-(1-p)^s)^{n-s}\geq
\label{expectation}
\end{equation}
$$
 e^{2\left(1-\frac{\gamma+2}{\gamma+1}\alpha\right)^2n^{\alpha}\ln^2 n(C_2-C_2^2)(1+o(1))-n^{1-2C_1\left(1-\frac{\gamma+2}{\gamma+1}\alpha\right)}(1+o(1))}\to\infty.
$$
It remains to prove that $\frac{{\sf Var}\tilde W(a)}{({\sf E}\tilde W(a))^2}\to 0$.\\

Consider distinct $s$-subsets $\mathcal{T}_d\subseteq V_n$, $d\in\{1,\ldots,{n\choose s}\}$, and events 

\noindent $B_d=$`the subgraph induced on $\mathcal{T}_d$ is isomorphic to $W_a^{\gamma}$', 

\noindent $\tilde B_d=B_d\,\wedge$ `every vertex outside $\mathcal{T}_d$ has a neighbor inside'.

Clearly, 
$$
{\sf E}\tilde W^2(a)={\sf E}\tilde W(a)+\sum_{d\neq\tilde d}{\sf P}(\tilde B_d\wedge \tilde B_{\tilde d})\leq
$$
\begin{equation}
{\sf E}\tilde W(a)+\sum_{\ell=0}^{s-1}\left[\left(1-(1-p)^s(2-(1-p)^{s-\ell})\right)^{n-2s+\ell}\sum_d{\sf P}(B_d)\sum_{\tilde d:\,|\mathcal{T}_d\cap\mathcal{T}_{\tilde d}|=\ell}{\sf P}(B_{\tilde d}|B_d)\right].
\label{variance}
\end{equation}

Fix $d\in\{1,\ldots,{n\choose s}\}$ and $\ell\in\{0,\ldots,s-1\}$. Let $\mathcal{S}_{\ell}\subset{V_n\choose{2s-\ell}}$ be the set of all $2s-\ell$-subsets of $V_n$ containing $\mathcal{T}_d$. For $S\in\mathcal{S}_{\ell}$, denote $\mathcal{D}[S]=:\mathcal{D}$ the set of all $s$-sets $\mathcal{T}_{\tilde d}$ such that $\mathcal{T}_d\cup\mathcal{T}_{\tilde d}=S$. Clearly,
$$
 \sum_{\tilde d:\,|\mathcal{T}_d\cap\mathcal{T}_{\tilde d}|=\ell}{\sf P}(B_{\tilde d}|B_d)=\sum_{S\in\mathcal{S}_{\ell}}\sum_{\tilde d:\,\mathcal{T}_{\tilde d}\in\mathcal{D}[S]}
 {\sf P}(B_{\tilde d}|B_d).
$$

Below, we estimate the value of $\sum_{\tilde d:\,\mathcal{T}_{\tilde d}\in\mathcal{D}}{\sf P}(B_{\tilde d}|B_d)$. The bound is given in Section 2.1.3.3. In Sections 2.1.3.1 and 2.1.3.2, we describe helpful construction and prove certain auxiliary bounds.\\

Fix $\epsilon>0$ as small as desired. Everywhere below, we distinguish two cases: {\it small $\ell$}, that is $\ell=o\left(n^{\alpha\frac{\gamma-9}{\gamma+1}}\right)$, and {\it large $\ell$}, that is $\ell\gg n^{\alpha\frac{\gamma-9}{\gamma+1}-\epsilon}$.

Let $G_1$ be the graph induced on $\mathcal{T}_d$.\\

First, consider $\ell=0$. In this case, clearly,  
$$
\sum_{\tilde d:\,\mathcal{T}_{\tilde d}\in\mathcal{D}}{\sf P}(B_{\tilde d}|B_d)=
\sum_{\tilde d:\,\mathcal{T}_{\tilde d}\in\mathcal{D}}{\sf P}(B_{\tilde d})=p^{E_s}(1-p)^{{s\choose 2}-E_s-{\ell\choose 2}}\frac{s!}{24^{[s-a]/[\gamma+1]-1}},
$$
where $E_s=\frac{\gamma+2}{\gamma+1}s-\frac{1}{\gamma+1}a-2$ is the number of edges in $W^{\gamma}_a$.\\

Now, let $\ell\in\{1,\ldots,s-1\}$.\\

Notice that $\sum_{\tilde d:\,\mathcal{T}_{\tilde d}\in\mathcal{D}}{\sf P}(B_{\tilde d}|B_d)$ equals summation over all possible ways of choosing an $\ell$-subset of $\mathcal{T}_d$ of `the number of ways of drawing remaining edges in the second $s$-set giving a copy of $W_a^{\gamma}$' multiplied by $p^{\text{number of remaining edges}}$. In Section 2.1.3.1, we describe possible ways of choosing an $\ell$-subset of $\mathcal{T}_d$. In Section 2.1.3.2, we  estimate from above both thee number of remaining edges and the number of ways of drawing them.\\

{\bf 2.1.3.1 Choosing an $\ell$-subset in $\mathcal{T}_d$}\\

\underline{If $\ell$ is small}, we ask the following question: What is the most likely structure of an $\ell$-vertex induced subgraph of $G_1$? More precisely, let us call a tree {\it $\gamma$-subdivided star}, if it is obtained from a star by subdividing every one of its edges by $\gamma$ vertices (i.e., every edge becomes $P_{\gamma+2}$). Note that $P_{\gamma+2}$ is a trivial $\gamma$-subdivided star. Given a graph $H$, let us call {\it the complexity of $H$} the maximum number of leaves in a vertex disjoint union of $\gamma$-subdivided substars in $H$.

Let us remind that we assume that $G_1$ is isomorphic to $W_a^{\gamma}$. Fix $x\in\mathbb{Z}_+$. Let us estimate the number of $\ell$-sets $\Upsilon\subset\mathcal{T}_d$ such that $G_1|_{\Upsilon}$ has a complexity at least $x$. If a $P_{\gamma+2}$ in $W_a^{\gamma}$ does not meet the $F_1$-part, then every its vertex in the graph has a degree at most $6$. Therefore, there are at most $s^6$ such paths. If a $P_{\gamma+2}$ meets the $F_1$-part, then it does not meet the $F_2$-part. Every such path consists of three segments (some of them may be trivial or empty): the first one and the last one do not meet the $F_1$-part, but the second one, in contrast, belongs to the $F_1$-part. There are at most $s^2$ such paths. Then, clearly, the desired number of $\ell$-sets is at most
\begin{equation}
 R_1(x):=\sum_{\tilde x=x}^{\lfloor\ell/(\gamma+1)\rfloor}s^{8\tilde x}{s-\tilde x(\gamma+1)\choose \ell-\tilde x(\gamma+1)}.
\label{small_1}
\end{equation}
Since 
$$
\left[s^{8(\tilde x+1)}{s-(\tilde x+1)(\gamma+1)\choose \ell-(\tilde x+1)(\gamma+1)}\right]/\left[s^{8\tilde x}{s-\tilde x(\gamma+1)\choose \ell-\tilde x(\gamma+1)}\right]\leq \frac{\ell^{\gamma+1}}{s^{\gamma-7}}=o(1),
$$
\begin{equation}
R_1(x)=s^{8x}{s-x(\gamma+1)\choose \ell-x(\gamma+1)}\left(1+O\left(\frac{\ell^{\gamma+1}}{s^{\gamma-7}}\right)\right).
\label{R_1}
\end{equation}

It means that almost all $\ell$-vertex (when $\ell$ is small enough) subgraphs of $W_a^{\gamma}$ does not have $\gamma$-subdivided stars. Therefore, given an $\ell$-vertex subgraph of $G_1$, we may assume (later, we will do it in a more formal way) that it is a forest. Moreover, in what follows, we show that we may even assume that it is an empty graph.\\

\underline{If $\ell$ is large}, then we do not need to look on subdivided stars --- in our computations, we just use the number of $\ell$-vertex subsets in $\mathcal{T}_d$, 
\begin{equation}
{s\choose\ell}.
\label{large_1}
\end{equation} 

{\bf 2.1.3.2 Choosing remaining edges}\\

Choose a set of $\ell$ vertices in $\mathcal{T}_d$ and denote it by $\Upsilon$ (at the moment, this set is arbitrary --- it may contain $\gamma$-subdivided stars). Denote by $\mathcal{G}$ and $G_2$ the graphs induced on $\Upsilon$ and $\mathcal{T}_{\tilde d}:=\Upsilon\cup[S\setminus \mathcal{T}_d]$ respectively. Let us estimate the probability ${\sf P}(B_{\tilde d}|B_d)$ that $G_2$ is isomorphic to $W_a^{\gamma}$. Roughly speaking, we need to compute the number of ways of embedding the vertices of $\mathcal{T}_{\tilde d}$ in $W_a^{\gamma}$ (put it differently, {\it building} $W_a^{\gamma}$ on the set $\mathcal{T}_{\tilde d}$) and multiply it by $p^{E}$, where $E$ is the number of extra edges we need for such an embedding. Since $E=2s-e(\mathcal{G})$, where $e(\mathcal{G})$ is the number of edges in $\mathcal{G}$, our computations rely heavily on an assumption about an edge-structure of $\mathcal{G}$.\\

Below, we use the following order $\preccurlyeq$ on the set of vertices of $W_a^{\gamma}$. On $V(F_1)$ and $V(F_2)$ it induces the tree order of the respective rooted trees (roots are the minimum vertices). Every $P_{\gamma+2}$ between the $F_1$- and $F_2$-parts is rooted in the vertex from $F_2$, and $\preccurlyeq$ also induced on the rooted path the tree order. For every $v\in V(F_1)$ and $u\notin V(F_1)$, $v\preccurlyeq u$. Let $v_1^0\preccurlyeq v_2^0\preccurlyeq \ldots\preccurlyeq v_a^0$ be the vertices of $F_1$, and $v_i^0\preccurlyeq v_i^1\preccurlyeq\ldots\preccurlyeq v_i^{\gamma+1}$ be the vertices of the $i$-th $P_{\gamma+2}$ (the order of the paths is arbitrary) between $F_1$ and $F_2$. We embed the vertices of $\mathcal{T}_{\tilde d}$ one by one in accordance to $\preccurlyeq$: at step $\tau$, we choose a vertex of $\mathcal{T}_{\tilde d}$ and map it to $v_{\tau}^0$, if $\tau\leq a$, and to $v_i^j$, where $i=1+\lfloor(\tau-a)/(\gamma+1)\rfloor$ and $j=\tau-a-(\gamma+1)(i-1)$, otherwise. Thus, we begin with $F_1$.\\

\begin{enumerate}

\item {\bf The $F_1$-part}


Here, we embed $a$ vertices of $\mathcal{T}_{\tilde d}$ in $F_1$. Recall that, on $V(F_1)$, $\preccurlyeq$ induces the linear tree order. Let us embed $i\in\{0,\ldots,\min\{a,\ell\}\}$ vertices inside $\Upsilon$ and $a-i$ vertices outside $\Upsilon$.\\

\underline{Let $\ell$ be small.}

First, let us assume that $i\geq 1$. All the vertices in $\Upsilon$ having {\it large} degrees (greater than $6$) should be embedded. Let $d_1,\ldots,d_y$ ($0\leq y\leq i$) be the degrees of these vertices. 
 The number of ways of choosing $i$ vertices (the vertices with large degrees should be among them) inducing $t$ disjoint paths in $\mathcal{G}$ is at most 
 \begin{equation}
 {i\choose t}\ell^t d_1\ldots d_y 6^{i-t-y},
\label{small_2}
\end{equation}
 and {\bf the number of edges in the union of these segments equals $i-t$}. Moreover, we should estimate the number of ways of embedding this union in $F_1\cong P_a$. It is at most 
\begin{equation}
 a^{t-1}(a-i+1)
\label{small_3}
\end{equation} 
($a^{t-1}$ ways of choosing the minimum vertices of the first $t-1$ paths and $a-i+1$ ways of choosing the minimum vertex of the last $t$-th path). There are 
\begin{equation}
\frac{(s-\ell)!}{(s-\ell-(a-i))!}
\label{small_4}
\end{equation} 
ways of choosing $a-i$ vertices to embed them in the remaining part of $F_1$. 
 
 Let us denote the constructed $P_a$ by $F^{\tilde d}_1$. Denote by $B$ and $A$ the subgraph of $\mathcal{G}$ induced on $\Upsilon\cap V(F^{\tilde d}_1)$ and the subgraph induced on $\Upsilon\setminus V(F^{\tilde d}_1)$ respectively. 

Clearly, $A$ is a forest (otherwise, our choice of $F^{\tilde d}_1$ was wrong). Let us remove from $A$ the trees that have neighbors in $V(B)$. Denote the remaining forest by $\tilde A$. Let $j$ be the number of tree components in $\tilde A$.
 
Let $x$ be the number of paths $P_{\gamma+2}$ in $\mathcal{G}$ such that the first vertex of every path is in $B$, and all the others are in $A$. (Notice that all the inner vertices of all these paths have degrees $2$ in $\mathcal{G}$ --- otherwise, our choice of $\mathcal{T}_{\tilde d}$ was wrong). Consider all the trees in $A$ that meet these paths, and remove the inner vertices of these paths. Denote the obtained forest by $A_0$. Let $j_0$ be the number of trees in $A_0$. Clearly, {\bf the number of edges in $\Upsilon$ equals exactly $\ell-j-j_0-t+x$}.\\

If $i=0$, then $\mathcal{G}$ is a disjoint union of $j$ trees. In this case, {\bf the number of edges in $\mathcal{G}$ equals $\ell-j$}. The number of ways of choosing $P_a$ outside $\Upsilon$ is at most 
\begin{equation}
\frac{(s-\ell)!}{(s-\ell-a)!}.
\label{small_0_1}
\end{equation}

\underline{For large $\ell$}, we use the following simple bound for the number of the desired embeddings: 
\begin{equation}
a^i{\ell\choose i}\frac{(s-\ell)!}{(s-\ell-(a-i))!}.
\label{large_2}
\end{equation} 

It remains to embed all the vertices of $\mathcal{T}_{\tilde d}\setminus F_1^{\tilde d}$ in $W_a^{\gamma}$. In accordance to the described order of embedding, we begin with the $P_{\gamma+2}$ between the roots of $F_1$ and $F_2$, and, on the way, we start building a new $P_{\gamma+2}$ as soon as the previous $P_{\gamma+2}$ is completed.\\


\item{\bf Paths $P_{\gamma+2}$}

As before, we, first, assume that \underline{$\ell$ is small}.

Here, we start from introducing  a certain classification of vertices from $V(A_0)\cup V(\tilde A)$. Then we define images of the vertices from $V(A)$ that lie in connected components having vertices in $B$. At the end, we embed recursively all the other vertices of $\mathcal{T}_{\tilde d}$.\\

\noindent {\bf Classification.}

In every tree from $A_0\sqcup\tilde A$ we choose a single vertex and call it {\it a root} (there are at most 
\begin{equation}
\left(\frac{\ell-i}{j+j_0}\right)^{j+j_0}
\label{small_5}
\end{equation}
 ways of choosing roots), see Fig.~\ref{fig2}.

\begin{figure}[h!]
\center{\includegraphics[scale=0.25]{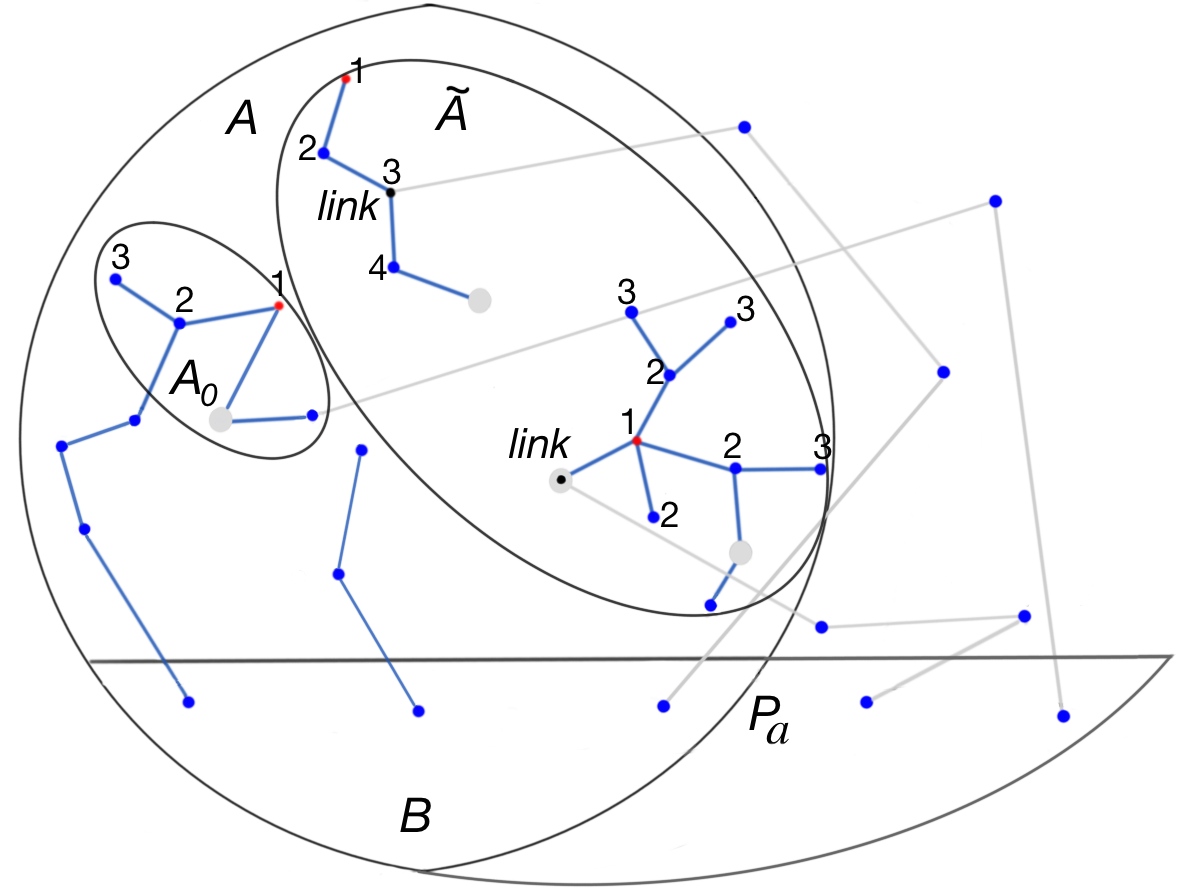}}
\caption{Labelling vertices of $A$: the links are colored black, the roots are colored red, the tails are colored grey; digits express the tree-orders. Grey edges appear within the described process; initial edges are colored blue.}
\label{fig2}
\end{figure}

Let us call a path of a rooted tree ${\it pendant}$, if the root does not belong to the set of the inner vertices of the path, its inner vertices have degree $2$ in the tree, the minimum vertex of the path (in accordance to the tree order of the given rooted tree) is either the root or has degree at least 3, and the maximum vertex has degree $1$. In every pendant path in $A_0\sqcup\tilde A$ choose at most one not minimum vertex, and call it {\it a tail} (there are at most 
\begin{equation}
2^{\ell-i-j-j_0}
\label{small_6}
\end{equation}
ways of choosing tails). Given a tree component $T^*$ (rooted in $r^*$) from $A_0\sqcup\tilde A$, define an order $\preccurlyeq_{T^*}$ on $V(T^*)$ in the following way. Let $\leq$ be the tree order of $T^*$. Let $v_1\leq v_2$ be two vertices of $T^*$. If there is no tail on the chain $r^*\leq\ldots\leq v_2$, then $v_1\preccurlyeq_{T^*}v_2$. If there is a tail $u$, and $u\leq v_1$, then $v_2\preccurlyeq_{T^*}v_1$. If $v_1<u\leq v_2$, and there is no vertex between $v_1$ and $u$, then $v_2\preccurlyeq_{T^*}v_1$ as well. Otherwise, $v_1$ and $v_2$ are incomparable. We would construct an embedding that preserves the tree orders (that is, $\preccurlyeq$ on $G_2\cong W_a^{\gamma}$ would induce $\preccurlyeq_{T^*}$ on $T^*$).

Finally, in every tree from $\tilde A$, we choose one more vertex, and call it {\it a link} (there are at most 
\begin{equation}
\left(\frac{\ell-i}{j}\right)^{j}
\label{small_7}
\end{equation} 
ways of choosing links). We make the following restriction on a choice of links: a link can not be a successor of a tail; if, in a tree $T^*$, there is a tail $u$ and a link $v$ {\it below} $v$ (i.e., $v\preccurlyeq_{T^*}u$), then $v$ is the only leaf below the tail in this tree.  
In what follows, we introduce an iterative algorithm of building the missing segments of paths $P_{\gamma+2}$ in $G_2$. A link in a tree $T^*$ from $\tilde A$ is the vertex that belongs to the first (in this iterative process) built edge (that belongs to a $P_{\gamma+2}$) between $T^*$ and $\overline{\Upsilon}$.\\

\noindent{\bf Building paths.}

As mentioned above, at every step $\tau\geq a+1$, we embed a vertex of $\mathcal{T}_{\tilde d}\setminus V(F_1^{\tilde d})$ in a {\it current} path $P_{\gamma+2}$. We proceed with a new path as soon as the current path is finished. All the minimum vertices of the paths are already embedded. On the way, we build paths of three following types: 1) paths in $\mathcal{G}$ connecting a vertex from $V(B)$ with a vertex from $V(A_0)$ (in fact, they are already built; for every such path, we should just define $\gamma+1$ successive steps devoted to its embedding --- there are at most 
\begin{equation}
s^x
\label{small_8}
\end{equation} 
ways of doing that, since there are $x$ such paths); 2) paths having non-trivial initial segments in $\mathcal{G}\setminus [\tilde A\sqcup A_0]$ with a minimum vertex lying $V(B)$; 3) all the other paths --- that do not have non-trivial initial segments in $\mathcal{G}$ (they may start either from a vertex in $V(B)$, or from a vertex in $V(F_{1}^{\tilde d})\setminus V(B)$). 

Let $d_0$ and $\Delta$ be the number of vertices in $V(A)$ having neighbors with {\it small} (at most $6$) degrees in $V(B)$ and the number of vertices in $V(A)$ having neighbors with large degrees in $V(B)$ respectively. Clearly, $d_1+\ldots+d_y-2y\leq \Delta\leq d_1+\ldots+d_y$. Let us estimate the number of ways of define successive steps devoted to embedding every {\it currently existed} segment (in $\mathcal{G}$) of a path of type 2). Clearly, there are at most 
\begin{equation}
{\lfloor s/(\gamma+1)\rfloor \choose d_0+\Delta}(d_0+\Delta)!
\label{small_9}
\end{equation}
 ways of doing that, since there are $d_0+\Delta$ such paths and less than $s/(\gamma+1)$ their potential final vertices.\\

It remains to built the remaining final segments of paths of type 2) and paths of type 3) completely.

Let us define the notions of active vertices, linked trees, and considered vertices. At the beginning of step $a+1$, the trees in $A_0$ are {\it linked}, the links and all the vertices from $\mathcal{T}_{\tilde d}\setminus[\Upsilon\cup V(F_1^{\tilde d})]$ are {\it active}, all the vertices from $V(A)$ connected in $\mathcal{G}$ by a path of {\it length} (i.e., the number of edges in the path) at most $\gamma+1$ with $V(B)$ and all the vertices from $V(F_1^{\tilde d})$ are {\it considered}. Note that, for every linked tree $T^*$, the order $\preccurlyeq_{T^*}$ gives a unique embedding of this tree in $W_a^{\gamma}$ (denote this embedding by $f_{T^*}$).

At step $\tau=a+1+\kappa(\gamma+1)$, $\kappa\in\mathbb{Z}_+$, we start to build the $(\kappa+1)$-th path $P\subset W_a^{\gamma}$. If this is a path of type 2), then we immediately move to the step $a+1+\kappa(\gamma+1)+q$, where $q$ is the number of edges in the respective segment in $\mathcal{G}$. We call the final vertex of this segment (which is in $V(A)$) the {\it current vertex}. If $P$ is a path of type 3) (i.e., $q=0$), then {\it the current vertex} is the pre-image of the minimum vertex of this path (it belongs to $V(F_{\tilde d}^1)$). 

Let $u\in V(F_2)$ be the maximum vertex of $P$. Two situations may happen: either there is a linked tree $T^*$ and a vertex $v\in V(T^*)$ such that $f_{T^*}(v)=u$, or not. If not, we map to $u$ an arbitrary active vertex $v$ and {\it deactivate} it. If this is the case, the active vertex may belong to a non-trivial tree $T^*$ from $\tilde A$. Independently of how $T^*$ appeared, find a tail (if there is one) $v_0$ such that either $v\preccurlyeq_{T^*} v_0$ or $v$ is the successor of $v_0$. We should build a $P_{\gamma+2-q-q_0}$ connecting the current vertex with $v^*$, where $v^*$ is the minimum vertex in $T^*$ below $v$, and $q_0$ is the number of vertices of $T^*$ but $v$ below $v$. 

If the deactivated vertex $v$ is either an isolated vertex from $\tilde A$ or a vertex outside $\Upsilon$, then it remains to build a $P_{\gamma+2-q}$ connecting the current vertex with $v$ (in this case, set $q_0=0$).

It remains to choose $\gamma+1-q-q_0$ vertices of $P$. We do it step by step, choosing an active vertex at every step and deactivating it. These active vertices should be either from $\mathcal{T}_{\tilde d}\setminus\Upsilon$, or end vertices of path components (not necessarily non-trivial) in $\tilde A$. In the latter case, the other end vertex of such a path should be its root, and the successive vertex should be its tail. Moreover, the number of vertices in this path should not be bigger than needed. In this latter case, we skip the number of moves equals to the number of edges in the path, and the root becomes {\it the current vertex}.

At the end of step $a+1+(\kappa+1)(\gamma+1)-1$, we {\it consider} all the vertices of the constructed path $P$. If $v$ belongs to a non-trivial tree $T^*$ from $\tilde A$, this tree becomes {\it linked} and we define (uniquely) its embedding $f_{T^*}$ in $W_a^{\gamma}$. Then we move to the step $a+1+(\kappa+1)(\gamma+1)$.

At the end of the very last step $s$, all the vertices of $G_2$ become considered, and all its vertices are embedded in $W_a^{\gamma}$.\\

Clearly, there are at most 
\begin{equation}
(s-\ell-(a-i)-d_0-\Delta+j)!
\label{small_10}
\end{equation} 
ways of making the construction.\\

\underline{For large $\ell$}, we use a similar (but much more simple) algorithm of the embedding. Here, we do not distinguish between trees in $A$, and choose a link, a root and tails in every tree in $A$. Assume that there are $j$ trees in $A$. Clearly, there are {\bf at most $\frac{\gamma+2}{\gamma+1}\ell-j$ edges} in $\mathcal{G}$. There are at most 
\begin{equation}
\left(\frac{\ell-i}{j}\right)^{2j}2^{\ell-i-j}
\label{large_3}
\end{equation} 
ways of choosing roots, links and tails, and at most 
\begin{equation}
(s-\ell-(a-i)+j)!
\label{large_4}
\end{equation} 
embeddings.

\end{enumerate}

{\bf 2.1.3.3 Estimating the variance}\\

We start from small $\ell$. In this case, we divide the summation into two parts w.r.t $i\geq 1$ and $i=0$. From~(\ref{small_1}),~(\ref{small_2})--(\ref{small_0_1}),~(\ref{small_5})--(\ref{small_10}), we get
$$
 \sum_{\tilde d:\,\mathcal{T}_{\tilde d}\in\mathcal{D}}
 {\sf P}(B_{\tilde d}|B_d)\leq p^{E_s}(1-p)^{{s\choose 2}-E_s-{\ell\choose 2}}\left(\sum_{x=0}^{\lfloor\ell/[\gamma+1]\rfloor}R_1(x)s^x \sum_{i=1}^{\min\{a,\ell\}}(a-i+1)\frac{(s-\ell)!}{(s-\ell-(a-i))!}\times\right.
$$
$$
\sum_{t=1}^i{i\choose t}\max_{y\in\{0,\ldots,i\},\,d_1,\ldots,d_y\in\{7,\ldots,\ell\}} \ell^t d_1\ldots d_y 6^{i-t-y}a^{t-1}\times
$$
$$
\max_{d_0\in\{0,\ldots,6i\},\,-2y\leq \Delta-d_1-\ldots-d_y\leq 0}{\lfloor s/(\gamma+1)\rfloor\choose d_0+\Delta}(d_0+\Delta)!\times
$$
$$
\max_{j,j_0:\,j+j_0\in\{1,\ldots,\ell-i\}}
\left(\frac{\ell-i}{j}\right)^j \left(\frac{\ell-i}{j+j_0}\right)^{j+j_0} 2^{\ell-i-j-j_0}\frac{(s-\ell-(a-i)-d_0-\Delta+j)!}{24^{([s-a]/[\gamma+1]-1-4(\ell-j-j_0))/4}}
\left(\frac{p}{1-p}\right)^{-[\ell-j-j_0-t+x]}
$$
$$
+\left.{s\choose\ell}\max_{j\in\{1,\ldots,\ell\}}\left(\frac{\ell}{j}\right)^{2j} 2^{\ell-j}\frac{(s-\ell+j)!}{24^{([s-a]/[\gamma+1]-1-4(\ell-j))/4}}
\left(\frac{p}{1-p}\right)^{-[\ell-j]}\right).
$$

Let 
$$
\xi_1(j,j_0,i)=\left(\frac{\ell-i}{j}\right)^j \left(\frac{\ell-i}{j+j_0}\right)^{j+j_0} 2^{-j-j_0}\frac{(s-\ell-(a-i)-d_0-\Delta+j)!}{24^{j+j_0}}\left(\frac{p}{1-p}\right)^{j+j_0}.
$$
Then, for $j+j_0\leq\ell-i-1$,
$$
 \frac{\xi_1(j,j_0+1,i)}{\xi_1(j,j_0,i)}=\frac{\ell-i}{j+j_0+1}\frac{1}{(1+1/(j+j_0))^{j+j_0}}\frac{1}{48}\frac{p}{1-p}<p\ell=o(1),\quad
 \frac{\xi_1(j+1,j_0,i)}{\xi_1(j,j_0,i)}=
$$
$$
 \frac{(\ell-i)^2}{(j+1)(j+j_0+1)}\frac{s-\ell-(a-i)-\Delta-d_0+j+1}{48(1+1/j)^j(1+1/(j+j_0))^{j+j_0}}\frac{p}{1-p}>
 \Theta(\ln n).
$$
Therefore, 
\begin{equation}
\max_{j,j_0:\,j+j_0\in\{1,\ldots,\ell-i\}}\xi_1(j,j_0,i)=\xi_1(\ell-i,0,i).
\label{max1}
\end{equation}
In particular,
\begin{equation}
\max_{j\in\{1,\ldots,\ell\}}\xi_1(j,0,0)=\xi_1(\ell,0,0)=2^{-\ell}\frac{(s-a-d_0-\Delta)!}{24^{\ell}}p^{\ell}(1+o(1)).
\label{max1.1}
\end{equation}

Let 
$$
\xi_2(d_0+\Delta)={\lfloor s/(\gamma+1)\rfloor\choose d_0+\Delta}(d_0+\Delta)!(s-a-d_0-\Delta)!.
$$
Then $\xi_2(d_0+\Delta+1)/\xi_2(d_0+\Delta)=\frac{s/(\gamma+1)-d_0-\Delta}{s-a-d_0-\Delta}\sim\frac{1}{\gamma+1}<1$. Therefore, 
\begin{equation}
\max_{d_0,\Delta}\xi_2(d_0+\Delta)=\xi_2(d_1+\ldots+d_y-2y).
\label{max2}
\end{equation}

Let 
$$
\xi_3(y;\,d_1,\ldots,d_y)=6^{-y}d_1\ldots d_y\xi_2(d_1+\ldots+d_y-2y).
$$
Then 
$$
\frac{\xi_3(y+1;\,d_1,\ldots,d_{y+1})}{\xi_3(y;\,d_1,\ldots,d_y)}\sim\frac{d_{y+1}}{6(\gamma+1)^{d_{y+1}-2}}<1
$$
since $d_{y+1}\geq 7$. Therefore, 
\begin{equation}
\max_{y,d_1,\ldots,d_y}\xi_3(d_1,\ldots,d_y)=1
\label{max3}
\end{equation} 
is achieved when $y=0$.

So, from~(\ref{max1})--(\ref{max3}),
$$
 \sum_{\tilde d:\,\mathcal{T}_{\tilde d}\in\mathcal{D}}
 {\sf P}(B_{\tilde d}|B_d)\leq p^{E_s}(1-p)^{{s\choose 2}-E_s-{\ell\choose 2}} \left(\sum_{x=0}^{\lfloor\ell/[\gamma+1]\rfloor}R_1(x)s^x
\sum_{i=1}^{\min\{a,\ell\}}(a-i+1)\right. \sum_{t=1}^i{i\choose t}\ell^t6^{i-t}a^{t-1}\times
$$
$$
\frac{(s-i)!}{24^{([s-a]/[\gamma+1]-1-4i)/4}}
\left(\frac{p}{1-p}\right)^{-[i-t+x]}
+\left.{s\choose\ell}\frac{s!}{24^{[s-a]/[\gamma+1]-1}}\right)(1+o(1)).
$$

Let $\eta_1(t)={i\choose t}\ell^t6^{-t}a^{t-1}\left(\frac{p}{1-p}\right)^{t}$. Then 
$$
\frac{\eta_1(t+1)}{\eta_1(t)}=\frac{i-t}{t+1}\ell \frac{a}{6}\frac{p}{1-p}\leq a^2\ell p=o(1).
$$
Therefore, $\sum_{t=1}^i\eta_1(t)=\eta_1(1)(1+o(1))$.

Let $\eta_2(i)=\eta_1(1)(a-i+1)6^i(s-i)!24^i\left(\frac{p}{1-p}\right)^{-i}$. Then
$$
 \frac{\eta_2(i+1)}{\eta_2(i)}=\frac{144(a-i)(1-p)}{(a-i+1)(s-i)p}=o(1).
$$
Therefore, $\sum_{i=1}^{\ell}\eta_2(i)=\eta_2(1)(1+o(1))$.

Let $\eta_3(x)=R_1(x)s^x\left(\frac{p}{1-p}\right)^{-x}$. Then, from~(\ref{R_1}), we get
$$
 \frac{\eta_3(x+1)}{\eta_3(x)}=\frac{(1-p)}{p}\frac{\ell^{\gamma+1}}{s^{\gamma-8}}\left(1+O\left(\frac{\ell^{\gamma+1}}{s^{\gamma-7}}\right)\right)=O\left(\frac{\ell^{\gamma+1}}{s^{\gamma-9}}\right),
$$
and the bound is uniform over all $x$. Therefore, $\sum_x \eta_3(x)=\eta_3(0)(1+o(1))$.

Putting it all together, we get
\begin{equation}
 \sum_{\tilde d:\,\mathcal{T}_{\tilde d}\in\mathcal{D}}
 {\sf P}(B_{\tilde d}|B_d)\leq p^{E_s}(1-p)^{{s\choose 2}-E_s-{\ell\choose 2}}{s\choose\ell}\frac{s!}{24^{[s-a]/[\gamma+1]-1}}(1+o(1)),
\label{conditional1}
\end{equation}
and the bound is uniform over all $d$ and $\mathcal{D}$.\\

If $\ell$ is large, then, similarly, from~(\ref{large_1}),~(\ref{large_2}),~(\ref{large_3}),~(\ref{large_4}),
$$
 \sum_{\tilde d:\,\mathcal{T}_{\tilde d}\in\mathcal{D}}
 {\sf P}(B_{\tilde d}|B_d)\leq{s\choose \ell} p^{E_s}(1-p)^{{s\choose 2}-E_s-{\ell\choose 2}}\left(\sum_{i=0}^{a}a^i{\ell\choose i}\right.\times
$$
$$
\left.\max_{j\in\{1,\ldots,\ell-i\}}\left(\frac{\ell-i}{j}\right)^{2j} 2^{\ell-i-j}\frac{(s-\ell+j)!}{24^{([s-a]/[\gamma+1]-1-4(\ell-j))/4}}
\left(\frac{p}{1-p}\right)^{-[\ell\frac{\gamma+2}{\gamma+1}-j]}\right)\leq
$$
\begin{equation}
{s\choose\ell}p^{E_s}(1-p)^{{s\choose 2}-E_s-{\ell\choose 2}}a^a{\ell\choose a}\frac{(s-a)!24^a}{24^{([s-a]/[\gamma+1]-1)/4}}
\left(\frac{p}{1-p}\right)^{-\ell\frac{1}{\gamma+1}-a}(1+o(1)).
\label{conditional2}
\end{equation}

From~(\ref{expectation}),~(\ref{variance}),~(\ref{conditional1}),~(\ref{conditional2}), we get
$$
\frac{{\sf E}\tilde W^2(a)}{({\sf E}\tilde W(a))^2}\leq
\frac{1}{{\sf E}\tilde W(a)}+
\frac{\sum_{\ell=0}^{\ell_0-1}{n-s\choose s-\ell}{s\choose\ell}\left(1-(1-p)^s(2-(1-p)^{s-\ell})\right)^{n-2s+\ell}(1-p)^{-{\ell\choose 2}}}{{n\choose s}(1-(1-p)^s)^{2(n-s)}}(1+o(1))+
$$
$$
 \frac{\sum_{\ell=\ell_0}^{s-1}{n-s\choose s-\ell}{s\choose\ell}\left(1-(1-p)^s(2-(1-p)^{s-\ell})\right)^{n-2s+\ell}(1-p)^{-{\ell\choose 2}}p^{-\ell\frac{1}{\gamma+1}-a}(a\ell)^a(s-a)!24^a}{{n\choose s}s!(1-(1-p)^s)^{2(n-s)}}(1+o(1)).
$$
Since $s(1-p)^s\to 0$ as $n\to\infty$, it is clear that
$$
 \frac{\left(1-(1-p)^s(2-(1-p)^{s-\ell}))\right)^{n-2s+\ell}}{(1-(1-p)^s)^{2(n-s)}}\sim e^{n(1-p)^{2s-\ell}},
$$
and the expression to the right approaches $0$ whenever $\ell=o(s)$. Below, we also use the inequality $\frac{{n-s\choose s-\ell}{s\choose\ell}}{{n\choose s}}\leq\left(\frac{es^2}{(n-s)\ell}\right)^{\ell}$. Thus, we get
$$
\frac{{\sf E}\tilde W^2(a)}{({\sf E}\tilde W(a))^2}\leq
\frac{1}{{\sf E}\tilde W(a)}+\frac{{n-s\choose s}}{{n\choose s}}(1+o(1))+\sum_{\ell=1}^{\ell_0-1}\left(\frac{es^2(1-p)^{-\ell/2}}{(n-s)\ell}\right)^{\ell}(1+o(1))+
$$
$$
\sum_{\ell=\ell_0}^{s-1}\left(\frac{es^2(1-p)^{-\ell/2}}{n\ell p^{1/[\gamma+1]}}(1+o(1))\right)^{\ell}e^{n(1-p)^{2s-\ell}}(1+o(1)).
$$
Notice that 1) ${n-s\choose s}\leq{n\choose s}$; 2) $(1-p)^{\ell/2}=1+o(1)$ whenever $\ell=o(n^{\alpha})$; 3) for $\ell\gg\frac{s}{\ln n}$, 
$$
 e^{n(1-p)^s}\leq e^{n^{1-2C_1(1-\frac{\gamma+2}{\gamma+1}\alpha)}}=e^{o(s)}=e^{o(\ell \ln n)};
$$
4) $\frac{(1-p)^{-\ell/2}}{\ell}$ is maximal when $\ell=s-1$; 5) $\frac{(1-p)^{-\ell/2}}{\ell}<1+o(1)$ whenever $\ell=o(s)$. These conclude the proof:
$$
\frac{{\sf E}\tilde W^2(a)}{({\sf E}\tilde W(a))^2}\leq 1+o(1)+\sum_{\ell=1}^{\ell_0-1}\left(n^{2\alpha-1+o(1)}\right)^{\ell}+\sum_{\ell=\ell_0}^{\lfloor\frac{s}{\ln\ln n}\rfloor}\left(n^{(2+1/[\gamma+1])\alpha-1+o(1)}\right)^{\ell}+
$$
$$
\sum_{\ell=\lfloor\frac{s}{\ln\ln n}\rfloor+1}^{s-1}\left(n^{(1+1/[\gamma+1])\alpha-1+o(1)}(1-p)^{-s/2}\right)^{\ell}\leq 
$$
$$
1+o(1)+\sum_{\ell=1}^{\infty}\left(n^{(2+1/[\gamma+1])\alpha-1+o(1)}\right)^{\ell}+\sum_{\ell=1}^{\infty}\left(n^{-(1-C_2)(1-\alpha\frac{\gamma+2}{\gamma+1})+o(1)}\right)^{\ell}=1+o(1).
$$

\subsection{Non-convergence for $\alpha\geq\frac{1}{2}$}
\label{after_one_half}

\subsubsection{The graph sequence}

We start from defining a graph sequence $W^{\gamma,*}_a$. Set $\omega_r(a)=\frac{r^a-1}{r-1}$.

Consider two graphs $F_1$ and $F_2$ with linear orders on their sets of vertices $\preccurlyeq_{F_1}$ and $\preccurlyeq_{F_2}$ respectively. Let $R_1$ and $R_2$ be minimum elements in $V(F_1)$ and $V(F_2)$ respectively. {\it The distance} between two vertices $u,v$ in a graph $F$ with a linear order $\preccurlyeq_F$ on $V(F)$ equals $d$, if the number of vertices in the maximum chain $u\preccurlyeq_F u_1\preccurlyeq_F\ldots\preccurlyeq_F v$ between $u$ and $v$ equals $d+1$. Let $\mathcal{E}$ be the set of all pairs $(u,v)$ where $u\in V(F_1)$, $v\in V(F_2)$ and $u$, $v$ are at the same distance from $R_1$, $R_2$ in $F_1$, $F_2$ respectively. Then the {\it ordered $\gamma$-product} of $F_1,F_2$ is the graph $F_1\cdot^{\gamma}_{\preccurlyeq_{F_1},\preccurlyeq_{F_2}} F_2$ obtained from the disjoint union $F_1\sqcup F_2$ in the following way: for every pair $(u,v)\in\mathcal{E}$, we add to the graph a $P_{\gamma+2}$ connecting $u$ and $v$. 

As in Section~\ref{before_one_half}, consider a simple path $F_1\cong P_a$ rooted at one of its end-points $R_1$, and a perfect $r$-ary tree $F_2$ having $\omega_r(a)$ vertices (its depth equals $a-1$) rooted at the only vertex having degree $r$. First, we consider the (not ordered) product $W_a^{\gamma}=F_1\cdot^{\gamma} F_2$. Let $\tilde F_1\cong P_{\omega_r(a)}$ be a simple path rooted at one of its end-points $\tilde R_1$, and $\tilde F_2$ be a perfect $r$-ary tree having $\omega_r(\omega_r(a))$ vertices (its depth equals $\omega_r(a)-1$) rooted at the only vertex $\tilde R_2$ having degree $r$. Second, consider the (not ordered) product $\tilde W_a^{\gamma}=\tilde F_1\cdot^{\gamma} \tilde F_2$. Notice that we consider pairwise disjoint sets $V(F_1),V(F_2),V(\tilde F_1),V(\tilde F_2)$. Let $\preccurlyeq_{\tilde F_1}$ be the tree order of the rooted tree $\tilde F_1$ on $V(\tilde F_1)$. Define the linear order $\preccurlyeq_{\tilde F_2}$ on $V(\tilde F_2)$ in the following way: for every $u,v\in V(\tilde F_2)$ such that the distance between $\tilde R_2$ and $u$ is less than the distance between $\tilde R_2$ and $v$, we set $u\preccurlyeq_{\tilde F_2}v$; on a set of vertices of $\tilde F_2$ at the same distance from $\tilde R_2$, the order $\preccurlyeq_{\tilde F_2}$ is defined arbitrarily.

\begin{figure}[h!]
\center{\includegraphics[scale=0.2]{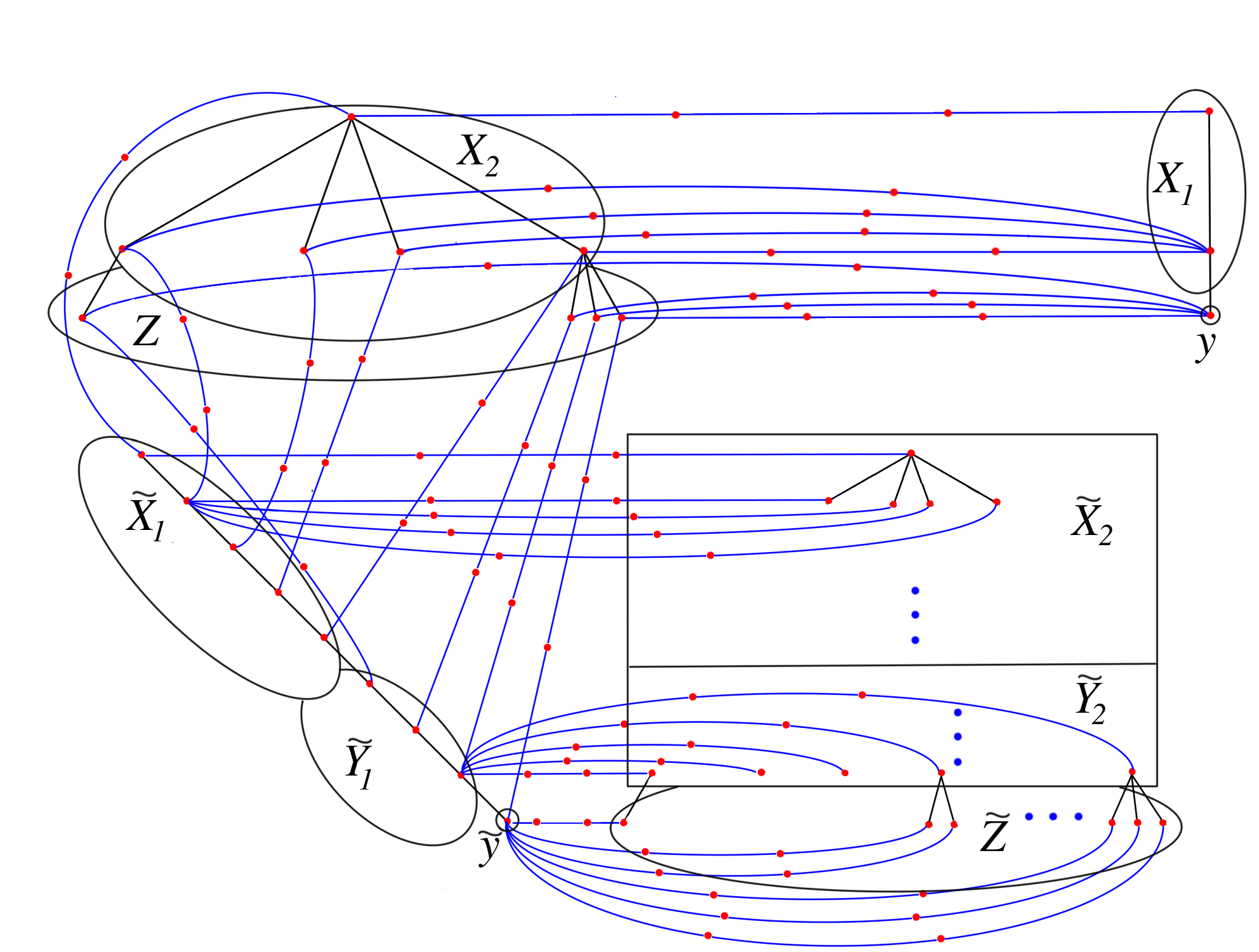}}
\caption{The graph $W^{\gamma,*}_a$ and the partition of its vertex set: $\gamma=2$, $r=4$, $a=2$, $i=21869$.}
\label{fig:upper_exponent}
\end{figure}

Finally, $W^{\gamma,*}_a$ is the union of $W_a^{\gamma}$, $\tilde W_a^{\gamma}$ and the ordered $\gamma$-product $\tilde F_1\cdot^{\gamma}_{\preccurlyeq_{\tilde F_1},\preccurlyeq_{\tilde F_2}} \tilde F_2$ (Fig.~\ref{fig:upper_exponent}). Clearly, $W^{\gamma,*}_a$ has 
$$
V_{\gamma,r}(a):=a+2(\gamma+1)\omega_r(a)+(\gamma+1)\omega_r(\omega_r(a))
$$ 
vertices and 
$$
E_{\gamma,r}(a):=a+2(\gamma+2)\omega_r(a)+(\gamma+2)\omega_r(\omega_r(a))-4
$$ 
edges.\\

\subsubsection{The graph process}

Then, we define a {\it $(\gamma,r)$ graph process} $W^{\gamma,*}(1)\subset W^{\gamma,*}(2)\subset\ldots$. In this process, there are arbitrarily large time moments when (arbitrarily large) graphs $W^{\gamma,*}_a$ appear. On the way between such two moments, we add vertices and edges to the previous $W^{\gamma,*}_a$ in order to obtain the next $W^{\gamma,*}_{a+1}$. We denote the vertices of the paths $F_1$ and $\tilde F_1$ of this $W^{\gamma,*}_a$ by $R_1=v_1\leq\ldots\leq v_a$ and $\tilde R_1=\tilde v_1\leq\ldots\leq\tilde v_{\omega_r(a)}$ respectively (here, $\leq$ are the respective tree orders). 

At time $i=1$, we have a graph $W^{\gamma,*}(0)\cong W^{\gamma,*}_1$ which is $P_{3\gamma+4}$ `joining' four roots $R_1,R_2,\tilde R_1,\tilde R_2$. 

At time $i$, we construct a graph $W^{\gamma,*}(i)$ on $(\gamma+1)(i+2)+a(i)$ vertices where $a(i)$ is the maximum integer such that $V_{\gamma,r}(a(i))\leq (\gamma+1)(i+2)+a(i)$. We assume that the graph $W^{\gamma,*}(i-1)$ is already constructed, and that is contains the above induced subgraph $W^{\gamma,*}_{a(i-1)}$ (with the roots $R_1,R_2,\tilde R_1,\tilde R_2$ and the vertex set $F_1\sqcup F_2\sqcup\tilde F_1\sqcup\tilde F_2$). Below, we set $a=a(i-1)$. The graph $W^{\gamma,*}(i)$ is obtained from $W^{\gamma,*}(i-1)$ in the following way.
\begin{itemize}
\item If 
$$
v\left(W^{\gamma,*}(i-1)\right)=V_{\gamma,r}(a),
$$
then $W^{\gamma,*}(i)$ is obtained from $W^{\gamma,*}(i-1)=W^{\gamma,*}_{a}$ by introducing one new vertex $v_{a+1}$ adjacent to the only vertex $v_a$ (the maximum vertex of $F_1$) of $W^{\gamma,*}_a$.

\item If, for some $j\in\{0,1,\ldots,r^{a}-1\}$, 
$$
v\left(W^{\gamma,*}(i-1)\right)=V_{\gamma,r}(a)+1+(2j+r^{\omega_r(a)}+\ldots+r^{\omega_r(a)+j-1})(\gamma+1),
$$
then $W^{\gamma,*}(i)$ is obtained from $W^{\gamma,*}(i-1)$ by introducing a neighbor $\tilde u$ of a leaf $u$ (if $a=1$, then $u=R_2$) of $F_2$ (such that, in $W^{\gamma,*}(i-1)$, the degree of $u$ is less than $3+r$, if $u\neq R_2$, and less than $2+r$, if $u=R_2$) and a $P_{\gamma+2}$ between $v_{a+1}$ and $\tilde u$.

\item If, for some $j\in\{0,1,\ldots,r^{a}-1\}$, 
$$
v\left(W^{\gamma,*}(i-1)\right)=V_{\gamma,r}(a)+1+(2j+r^{\omega_r(a)}+\ldots+r^{\omega_r(a)+j-1}+1)(\gamma+1),
$$
then $W^{\gamma,*}(i)$ is obtained from $W^{\gamma,*}(i-1)$ by introducing a neighbor $\tilde v_{\omega_r(a)+j+1}$ of the maximum vertex $\tilde v_{\omega_r(a)+j}$ of the simple path of $W^{\gamma,*}(i-1)$ extending $\tilde F_1$ (and this vertex becomes the new maximum vertex), and a $P_{\gamma+2}$ between $\tilde v_{\omega_r(a)+j}$ and the neighbor of the leaf of $F_2$ that was attached on the previous step.

\item If, for some $j\in\{0,1,\ldots,r^{a}-1\}$, $j_0\in\{0,1,\ldots,r^{\omega_r(a)+j}-1\}$, 
$$
v\left(W^{\gamma,*}(i-1)\right)=V_{\gamma,r}(a)+1+(2j+r^{\omega_r(a)}+\ldots+r^{\omega_r(a)+j-1}+2+j_0)(\gamma+1),
$$
then $W^{\gamma,*}(i)$ is obtained from $W^{\gamma,*}(i-1)$ by introducing a neighbor $\tilde u$ of a leaf $u$ (if $a=1$ and $j=0$, then $u=\tilde R_2$) of the perfect $r$-ary tree extending $\tilde F_2$ constructed at the moment $i-1-(j_0+2)$ (such that, in $W^{\gamma,*}(i-1)$, the degree of $u$ is less than $2+r$, if $u\neq \tilde R_2$, and less than $1+r$, if $u=\tilde R_2$) and a $P_{\gamma+2}$ between $\tilde v_{\omega_r(a)+j+1}$ and $\tilde u$.
\end{itemize}
We will call $a(i)$ {\it the floor} of $W^{\gamma,*}(i)$.\\

For appropriate $\gamma$ and $r$, the desired EMSO sentence expresses the $(\gamma,r)$-property defined below. We say that a graph $\mathcal{G}$ has the {\it $(\gamma,r)$-property}, if
\begin{center}
There exist an {\it even} positive integer $a$, a positive integer $i$ and a $(\gamma,r)$ process $W^{\gamma,*}(1)\subset W^{\gamma,*}(2)\subset\ldots\subset W^{\gamma,*}(i)$ such that:

\vspace{0.2cm}

$\bullet$ $W^{\gamma,*}(1)\subset\ldots\subset W^{\gamma,*}(i)\subset\mathcal{G}$ are induced subgraphs in $\mathcal{G}$;

\vspace{0.1cm}

$\bullet$ $a$ is the floor of $W^{\gamma,*}(i)$;

\vspace{0.1cm}

$\bullet$ there is no induced $W^{\gamma,*}(i+1)\subset\mathcal{G}$ such that $W^{\gamma,*}(1)\subset\ldots\subset W^{\gamma,*}(i+1)$ is a $(\gamma,r)$ process.
\end{center}

\subsubsection{The sentence}

Let us construct an EMSO sentence $\varphi$ that expresses the $(\gamma,r)$-property. Clearly, such a sentence can be written in the following way (see Fig.~\ref{fig:upper_exponent}):
$$
 \varphi=\exists X_1\exists X_2\exists\tilde X_1\exists \tilde X_2\exists \tilde Y_1\exists \tilde Y_2\exists Z\exists\tilde Z\exists\Gamma_1\ldots\exists\Gamma_7\exists y\exists\tilde y
$$
$$
 [\phi_{\gamma,r}(X_1,X_2,\Gamma_1)\wedge\phi_{\gamma,r}(\tilde X_1,\tilde X_2,\Gamma_3)\wedge\phi_{\gamma,r}(\tilde X_1\vee\tilde Y_1,\tilde X_2\vee\tilde Y_2,\Gamma_3\vee\Gamma_6)\wedge
 \phi^*(X_2,\tilde X_1,\Gamma_2)\wedge\phi^*(Z,\tilde Y_1,\Gamma_5)\wedge
$$
$$
  \mathrm{PATHS}(X_1,X_2,\tilde X_1,\Gamma_1,\Gamma_2)\wedge\mathrm{LAST}(X_1,Z,y,\Gamma_4)\wedge\mathrm{LAST}(\tilde Y_1,\tilde Z,\tilde y,\Gamma_7)\wedge
$$  
$$  
  \mathrm{LEAVES}(X_2,Z)\wedge\mathrm{LEAVES}(\tilde Y_2,\tilde Z)\wedge\mathrm{EVEN}(X_1)\wedge \mathrm{MAX}\wedge\varnothing\wedge\mathrm{EDGES}].
$$

Let $\phi_{\gamma,r}(X_1,X_2,\Gamma)$ be a FO formula saying that the graph induced on $[X_1]\cup[X_2]\cup[\Gamma]$ is isomorphic to $W_a^{\gamma}$ for some $a$ (for $\gamma=0$ and $r=4$, see the construction in Section~\ref{sentence_desc}~and~\cite{Zhuk_Le_Bars}, proof of Theorem 3; for arbitrary $\gamma$ and $r$, is admits a straightforward generalization). In particular, $X_1\cong F_1$, $X_2\cong F_2$, and the inner vertices of $\gamma+2$-paths belong to $\Gamma$.

Let $\phi^*(X_1,X_2,\Gamma)$ be a FO formula saying that the vertices of $[X_1]\cup[X_2]\cup[\Gamma]$ and edges of $[\Gamma]$ are covered by a disjoint union of $P_{\gamma+2}$ having ends in $[X_1]$ and $[X_2]$ (one end in one set) and inner vertices in $[\Gamma]$. The existence of such a formula is straightforward: one may say that, 1) for every vertex $x_1$ from $[X_1]$, there is the only vertex $x_2$ in $[X_2]$ which is connected with $x_1$ by $P_{\gamma+2}$ having all inner vertices in $[\Gamma]$, and vice versa; 2) every vertex in $[\Gamma]$ belongs to a component $P_{\gamma}$ having ends $y_1,y_2$ such that $y_1$ and $y_2$  have the only neighbors in $[X_1]\cup[X_2]$, and these neighbors are from different sets.

Let a FO formula $\mathrm{PATHS}(X_1,X_2,\tilde X_1,\Gamma_1,\Gamma_2)$ say that every inclusion-maximum set of vertices in $[X_2]$ mapped by the paths $P_{\gamma+2}$ having inner vertices in $[\Gamma_1]$ into one common vertex from $[X_1]$ has the following property. The images of these vertices in $[\tilde X_1]$ under the bijection produced by the paths $P_{\gamma+2}$ having inner vertices in $[\Gamma_2]$ induces a simple path in $[\tilde X_1]$. It is clear that such a formula exists under the assumption that $[\tilde X_1]$ is a simple path (one may say that there are only two vertices with degree 1, and all the others have degree 2 in the induced subgraph).

Let a FO formula $\mathrm{LAST}(X,Z,y,\Gamma)$ say that the vertex $y$ adjacent to the maximum (the root is defined by the respective formula $\phi_{\gamma,r}$) vertex of the path induced by $[X]$, and is not adjacent to any other vertex of $[X]$; $[\Gamma]$ induces a disjoint union of $P_{\gamma}$ having first and last vertices such all the first vertices (and only they) are adjacent to $y$, and every last vertex has exactly one neighbor in $[Z]$; every vertex from $[Z]$ has exactly one neighbor in $[\Gamma]$, and this neighbor is a last vertex.

A FO formula $\mathrm{LEAVES}(X,Z)$ says that every vertex in $[Z]$ has exactly one neighbor in $[X]$, this neighbor has a degree at most $1$ in the graph induced by $[X]$, every vertex in $[X]$ has at most $r$ neighbors in $[Z]$, and there are no edges in $[Z]$.

An EMSO formula $\mathrm{EVEN}(X)$ says that, under the condition that $[X]$ induces a path, the cardinality of $[X]$ is even. One monadic variable is enough to say this: one may say that, for some $X_0\subset[X]$, all edges in the subgraph induced by $[X]$ are between $X_0$ and $[X]\setminus X_0$, and the ends of the path belong to different sets.

A FO formula $\mathrm{MAX}$ with unary predicates $X_1,X_2,\tilde X_1,\tilde X_2,\tilde Y_1,\tilde Y_2,Z,\tilde Z,\Gamma_1,\ldots,\Gamma_7$ says the following:
\begin{itemize}
\item If every vertex having a degree at most $1$ in the graph induced on $[X_2]$ has $r$ neighbors in $[Z]$, and every vertex having a degree at most $1$ in the graph induced on $[\tilde Y_2]$ has $r$ neighbors in $[\tilde Z]$, then there is no neighbor of $y$ outside $U:=[X_1]\sqcup\ldots\sqcup[\Gamma_7]\sqcup\{y,\tilde y\}$ that has no other neighbors in $U$.

\item If there exists a vertex having a degree at most $1$ in the graph induced on $[\tilde Y_2]$ with at most $r-1$ neighbors in $[\tilde Z]$, then there is no vertex $w$ outside $U$ such that the following property holds: $w$ is adjacent to exactly one vertex $u$ from $U$, this neighbor $u$ belongs to the set of vertices having a degree at most $1$ in the graph induced on $[\tilde Y_2]$, $u$ has at most $r-1$ neighbors in $[Z]$, and there exists a path $P_{\gamma+2}$ connecting $w$ with $\tilde y$ and having inner vertices outside $U$ and not adjacent to any vertex from $U$ (the only exception is the vertex after $\tilde y$ which is adjacent only to $\tilde y$).

\item Finally, if every vertex having a degree at most $1$ in the graph induced on $[\tilde Y_2]$ has $r$ neighbors in $[\tilde Z]$, but there exists a vertex having a degree at most $1$ in the graph induced on $[X_2]$ with at most $r-1$ neighbors in $[Z]$, then there is no vertex $w$ outside $U$ such that the following property holds: $w$ is adjacent to exactly one vertex $u$ from $U$, this neighbor $u$ belongs to the set of vertices having a degree at most $1$ in the graph induced on $[X_2]$, $u$ has at most $r-1$ neighbors in $[Z]$, and there exists a path $P_{\gamma+2}$ connecting $w$ with $y$ and having inner vertices outside $U$ and not adjacent to any vertex from $U$ (the only exception is the vertex after $y$ which is adjacent only to $y$).
\end{itemize}

$\varnothing$ says that the sets $[X_1],[X_2],[\tilde X_1],[\tilde X_2],[\tilde Y_1],[\tilde Y_2],[Z],[\tilde Z],[\Gamma_1],\ldots,[\Gamma_7]$ are pairwise disjoint.

Finally, $\mathrm{EDGES}$ says that, for every unconsidered pair of sets, there are no edges between the vertices in these sets (e.g., there are no edges between $X_1$ and $\tilde X_1$).\\

Clearly, the second statement of Theorem~\ref{main} follows from the lemma below.

\begin{lemma}

Let $\alpha\in(0,1)$.

\begin{enumerate}

\item Let $\varepsilon>0$. Then a.a.s., for every $a$ such that $|V_{\gamma,r}(a)|\geq 2\left(1-\frac{\gamma+2}{\gamma+1}\alpha\right)n^{\alpha}\ln n+\varepsilon$, in $G(n,n^{-\alpha})$, there is no induced copy of $W_a^{\gamma,*}$.

\item Let $0<\beta<\min\left\{\alpha,\frac{2}{3}(1-\alpha)\right\}$. Then, for all large enough positive integer $\gamma$,
\begin{enumerate}
\item a.a.s., for every $a$ such that  $|V_{\gamma,r}(a)|\leq n^{\beta}$, in $G(n,n^{-\alpha})$, there exists an induced copy of $W_a^{\gamma,*}$;
\item a.a.s., for every set $V\subset V_n$ with $|V|\leq n^{\beta}$ and every pair of vertices $u,v\in V$, there exists an induced $P_{\gamma+1}=w_1\ldots w_{\gamma+1}$ outside $V$ such that the only neighbor of $w_1$ in $V$ is $u$, the only neighbor of $w_{\gamma+1}$ in $V$ equals $v$, vertices $w_2,\ldots,w_{\gamma}$ do not have neighbors in $V$.
\end{enumerate}

\end{enumerate}
\label{Lem2}
\end{lemma}

Indeed, let $0<\beta<\min\left\{\alpha,\frac{2}{3}(1-\alpha)\right\}$. To prove the result it is enough to show that, for large enough $r$, there are two sequences $n_i,\,i\in\mathbb{N},$ and $m_i,\,i\in\mathbb{N}$, with the following property. For all large enough $i$, there are an even number $a_1(i)$ and an odd number $a_2(i)$ such that 
$$
|V_{\gamma,r}(a_1(i))|\leq n_i^{\beta},\quad
|V_{\gamma,r}(a_1(i)+1)|>2\left(1-\frac{\gamma+2}{\gamma+1}\alpha\right)n_i^{\alpha}\ln n_i+\varepsilon;
$$
$$
|V_{\gamma,r}(a_2(i))|\leq m_i^{\beta},\quad
|V_{\gamma,r}(a_2(i)+1)|>2\left(1-\frac{\gamma+2}{\gamma+1}\alpha\right)m_i^{\alpha}\ln m_i+\varepsilon.
$$
It is clear that $n_i=2\left\lfloor\left[(\gamma+1)\frac{r^{\frac{r^{2i}-1}{r-1}}-1}{r-1}\right]^{1/\beta}\right\rfloor$,
$m_i=2\left\lfloor\left[(\gamma+1)\frac{r^{\frac{r^{2i+1}-1}{r-1}}-1}{r-1}\right]^{1/\beta}\right\rfloor$ are appropriate.

\subsubsection{Proof of Lemma~\ref{Lem2}}

Let $W^*(a)$ be the number of induced copies of $W^{\gamma,*}_a$ in $G(n,p)$. Let $s=V_{\gamma,r}(a)$. Clearly $E_{\gamma,r}(a)=\frac{\gamma+2}{\gamma+1}s-\frac{1}{\gamma+1}a-4$.

\noindent 1. Let $\varepsilon>0$.  For $s\geq 2(1-\frac{2+\gamma}{1+\gamma}\alpha)n^{\alpha}\ln n+\varepsilon$,
$$
 {\sf E}W^*(a)={n\choose s}\frac{s!}{(r!)^{([s-a-2(\gamma+1)\omega_r(a)]/[\gamma+1]-1)/r}}p^{\frac{\gamma+2}{\gamma+1}s-\frac{1}{\gamma+1}a-4}(1-p)^{{s\choose 2}-\frac{\gamma+2}{\gamma+1}s+\frac{1}{\gamma+1}a+4}\leq
$$
$$
 e^{s\left[\left(1-\frac{\gamma+2}{\gamma+1}\alpha\right)\ln n-\frac{sp}{2}-\frac{\ln(r!)}{r(\gamma+1)}+O(p)\right]+O(\ln n\ln\ln n)}.
$$
Therefore, ${\sf P}(W^*(a)>0)\leq{\sf E}W^*(a)\to 0$ as $n\to\infty$.

\medskip

\noindent 2. Now, let $0<\beta<\frac{1}{2}$. Then, for $\ln\ln n\ll s\leq n^{\beta}$,
$$
 {\sf E}W^*(a)\geq e^{s\left[\left(1-\frac{\gamma+2}{\gamma+1}\alpha\right)\ln n-\frac{\ln(r!)}{r(\gamma+1)}+o(1)\right]+O(\ln n\ln\ln n)}\to\infty.
$$
To prove 2.(a), it remains to show that, for such $s$, $\frac{{\sf Var}W^*(a)}{({\sf E}W^*(a))^2}\to 0$.\\

Consider distinct $s$-subsets $\mathcal{T}_d\subseteq V_n$, $d\in\{1,\ldots,{n\choose s}\}$, and events $B_d=$`the subgraph induced on $\mathcal{T}_d$ is isomorphic to $W_a^{\gamma,*}$'. Then 
$$
{\sf Var}W^*(a)={\sf E}W^*(a)+\sum_{\ell=0}^{s-1}\sum_{d\neq\tilde d:\,|\mathcal{T}_d\cap\mathcal{T}_{\tilde d}|=\ell}{\sf P}(B_d\wedge B_{\tilde d})-({\sf E}W^*(a))^2.
$$
Let $\ell\in\{1,\ldots,s-1\}$ and $\Upsilon$ be an $\ell$-subset of $W_a^{\gamma,*}$. Let $\mathcal{P}$ be the set of all $P_{\gamma+2}$ between $F_1, F_2$; $F_2,\tilde F_1$ and $\tilde F_1,\tilde F_2$. Let $x$ be the number of paths from $\mathcal{P}$ that are entirely in $W_a^{\gamma,*}|_{\Upsilon}$. Clearly, the number of edges in $W_a^{\gamma,*}|_{\Upsilon}$ is at most $\ell+x\leq \frac{\gamma+2}{\gamma+1}\ell$ since paths from $P\in\mathcal{P}$ can be ordered in a way $P_1,P_2\ldots$ such that, for every $i$, $|V(P_i)\cap [V(P_1)\cup\ldots\cup V(P_{i-1})]|\leq 1$. Note that $\frac{(1-p)^{-\ell/2}}{1-p}\leq 1+o(1)$ since the function to the left decreases in $\ell$ on $[1,s]$. From this,
$$
{\sf Var}W^*(a)\leq{\sf E}W^*(a)+\sum_{\ell=1}^{s-1}\sum_{d\neq\tilde d:\,|\mathcal{T}_d\cap\mathcal{T}_{\tilde d}|=\ell}{\sf P}(B_d\wedge B_{\tilde d})\leq
{\sf E}W^*(a)+
$$
$$
\sum_{\ell=1}^{s-1}{n\choose s}{s\choose\ell}{n-s\choose s-\ell}\frac{[s!]^2}{(r!)^{2([s-a-2(\gamma+1)\omega_r(a)]/[\gamma+1]-1)/r-\ell}} p^{2E_{\gamma,r}(a)-\frac{\gamma+2}{\gamma+1}\ell} (1-p)^{2\left({s\choose 2}-E_{\gamma,r}(a)\right)-{\ell\choose 2}+\frac{\gamma+2}{\gamma+1}\ell}
$$
$$
\leq{\sf E}W^*(a)+\sum_{\ell=1}^{s-1}\left[{\sf E}W^*(a)\right]^2\left(\frac{es^2}{(n-s)\ell}\right)^{\ell}(r!)^{\ell}p^{-\frac{\gamma+2}{\gamma+1}\ell} (1-p)^{-{\ell\choose 2}+\frac{\gamma+2}{\gamma+1}\ell}\leq
$$
$$
{\sf E}W^*(a)+\sum_{\ell=1}^{s-1}\left[{\sf E}W^*(a)\right]^2\left(\frac{r^r s^2(1-p)^{-\ell/2}}{n \ell p^{\frac{\gamma+2}{\gamma+1}}}(1+o(1))\right)^{\ell}\leq
$$
$$
{\sf E}W^*(a)+\sum_{\ell=1}^{s-1}\left[{\sf E}W^*(a)\right]^2\left(n^{-\frac{1}{3}(1-\alpha)+\frac{\alpha}{\gamma+1}+o(1)}\right)^{\ell}=o(({\sf E}W^*(a))^2)
$$
whenever $\gamma+1>\frac{3\alpha}{1-\alpha}$.\\

It remains to prove 2.(b). Let $\frac{\gamma+1}{\gamma+3}>\alpha$. Then, for some constant $C_1>0$, a.a.s., for every pair of vertices $u,v\in V_n$, there are at least $C_1 n^{\gamma+1}p^{\gamma+2}$ disjoint induced $P_{\gamma+3}$ in $G(n,n^{-\alpha})$ connecting $u$ and $v$ (see \cite{Spencer_Counting}, Theorem 2). Moreover, there exists $C_2>0$ such that, a.a.s., for every three vertices $u,v,w\in V_n$ there are at most $C_2 n^{\gamma+1}p^{\gamma+3}$ induced $P_{\gamma+3}$ in $G(n,n^{-\alpha})$ connecting $u$ and $v$ and having at least one neighbor of $w$ (see \cite{Spencer_Counting}, Theorem 2).

Clearly, if 
\begin{equation}
n^{\beta} \ll n^{\gamma+1}p^{\gamma+2},\quad n^{\gamma+1}p^{\gamma+3}n^{\beta}\ll n^{\gamma+1} p^{\gamma+2},
\label{small_set}
\end{equation}
then a.a.s., for every $V\subset V_n$ having $|V|=\lfloor n^{\beta}\rfloor$ and every pair of vertices $u,v\in V$, there exists an induced $P_{\gamma+2}=w_1\ldots w_{\gamma+1}$ outside $V$ such that the only neighbor of $w_1$ in $V$ is $u$, the only neighbor of $w_{\gamma+2}$ in $V$ equals $v$, vertices $w_2,\ldots,w_{\gamma+1}$ do not have neighbors in $V$. But the second condition in (\ref{small_set}) holds since $\beta<\alpha$, and the first one holds for all $\gamma>\frac{\beta+2\alpha-1}{1-\alpha}$. Lemma is proven.

\section*{Acknowledgements}

This work is supported by the grant N 18-71-00069 of Russian Science Foundation.

\renewcommand{\refname}{References}

\end{document}